\documentclass[a4paper,12pt]{article}

\usepackage{amsmath}
\usepackage{amsthm}
\usepackage{amssymb}  % defines \nmid, for example
\usepackage{latexsym} % for \Box
\usepackage[all]{xy}
\usepackage{comment}
\usepackage{rotating}
\usepackage{slashbox}
\usepackage{multirow}
\usepackage{rotating}

\newcommand{\Fp}{{\mathbb{F}_p}}

\newcommand{\Char}{\operatorname{char}}

\newcommand{\Q}{{\mathbb Q}}

\newfont{\wncyr}{wncyr10 at 12pt}
\newfont{\wncyrten}{wncyr10 at 10pt}

\newenvironment{Proof}{\par\noindent{\sc Proof:}}%
                      {\hspace*{\fill}\nobreak$\Box$\par\medskip}
                       {\hspace*{\fill}\nobreak$\Box$\par\medskip}
\newenvironment{myitemize}
{\begin{itemize}
\setlength{\itemsep}{1pt}
\setlength{\parskip}{0pt}
\setlength{\parsep}{0pt}}
{\end{itemize}}

\newtheorem{Proposition}{Proposition}[section]
\newtheorem{Theorem}[Proposition]{Theorem}
\newtheorem{Lemma}[Proposition]{Lemma}
\newtheorem{Corollary}[Proposition]{Corollary}

\theoremstyle{definition}

\newcounter{nootje}
\begin{document}
\normalsize
\title{Edwards Curves and Gaussian Hypergeometric Series}
\author{Mohammad Sadek and Nermine El-Sissi}
\date{}
\maketitle
\let\thefootnote\relax\footnote{Mathematics Subject Classification: 11T24, 14H52}
\begin{abstract}{\footnotesize Let $E$ be an elliptic curve described by either an Edwards model or a twisted Edwards model over $\Fp$, namely, $E$ is defined by one of the following equations $x^2+y^2=a^2(1+x^2y^2),\, a^5-a\not\equiv 0$ mod $p$, or, $ax^2+y^2=1+dx^2y^2,\,ad(a-d)\not\equiv0$ mod $p$, respectively. We express the number of rational points of $E$ over $\Fp$ using the Gaussian hypergeometric series $\displaystyle {_2F_1}\left(\begin{matrix}
                          \phi&\phi \\
                          {} & \epsilon
                        \end{matrix}\Big| x\right)$ where $\epsilon$ and $\phi$ are the trivial and quadratic characters over $\Fp$ respectively. This enables us to evaluate $|E(\Fp)|$ for some elliptic curves $E$, and prove the existence of isogenies between $E$ and Legendre elliptic curves over $\Fp$.}
\end{abstract}

\section{Introduction}
In \cite{Greene} Greene initiated the study of Gaussian hypergeometric series over finite fields. These series are analogous to the classical hypergeometric series. Several authors managed to find congruence relations satisfied by special values of these series, see \cite{Osburn}. Many special values of these series were determined.

One of the striking aspects of hypergeometric series is that some of their special values are linked to the number of rational points on some families of algebraic curves over finite fields. Two families of elliptic curves were discussed in \cite{OnoGauss}. The number of rational points of an elliptic curve $E$ described by a Legendre model, namely, $y^2=x(x-1)(x-\lambda),\,\lambda(\lambda-1)\not\equiv 0$ mod $p$, over the finite field $\Fp$ satisfies the following identity
 \[|E(\Fp)|=1+p+p\phi(-1)\cdot{_2F_1}\left(\begin{matrix}
                          \phi&\phi \\
                          {} & \epsilon
                        \end{matrix}\Big| \lambda\right)\]
                        where $\epsilon$ and $\phi$ are the trivial and quadratic characters over $\Fp$ respectively. The latter identity was used to evaluate the hypergeometric series $\displaystyle{_2F_1}\left(\begin{matrix}
                          \phi&\phi \\
                          {} & \epsilon
                        \end{matrix}\Big| \lambda\right)$ when $\lambda\in\{-1,1/2,2\}$. If $E$ is defined by a Clausen model, $y^2=(x-1)(x^2+\lambda),\lambda(\lambda+1)\not\equiv 0$ mod $p$, then \[\left(1+p-|E(\Fp)|\right)^2=p+p^2\phi(\lambda+1)\cdot{_3F_2}\left(\begin{matrix}
                          \phi&\phi &\phi\\
                          {} & \epsilon&\epsilon
                        \end{matrix}\Big|\frac{\lambda}{\lambda+1}\right).\]
Again this was exploited in order to evaluate the hypergeometric series $\displaystyle{_3F_2}\left(\begin{matrix}
                          \phi&\phi &\phi\\
                          {} & \epsilon&\epsilon
                        \end{matrix}\Big|\frac{\lambda}{\lambda+1}\right) $ at some specific values of $\lambda$.

More Gaussian hypergeometric series appear in formulas describing the number of rational points on higher genus curves over finite fields. Some of these formulas can be found in \cite{Barmanalgebraiccurves} where the following family of algebraic curves are discussed
\[y^l=x(x-1)(x-\lambda),\;\lambda(\lambda-1)\not\equiv 0\textrm{ mod }p,\;l\ge2.\]

In this note we are interested in elliptic curves described by Edwards models or twisted Edwards models, namely
\begin{eqnarray*}x^2+y^2&=&a^2(1+x^2y^2),\,a^5-a\not\equiv 0\textrm{ mod }p;\\
ax^2+y^2&=&1+dx^2y^2,\,ad(a-d)\not\equiv0\textrm{ mod }p
\end{eqnarray*}
respectively.

 Edwards models of elliptic curves were proposed in \cite{Edwards} and have been used since then in many cryptographic applications. The main advantage enjoyed by these curves is that the group law is simpler to state than on other models representing elliptic curves. In addition, any elliptic curve defined over an algebraically closed field $k$ can be expressed in the form $x^2+y^2=a^2(1+x^2y^2)$.
Twisted Edwards models appeared for the first time in \cite{Bernstein} to express more elliptic curves over finite fields with the addition law being easily formulated.

 We show that the number of rational points on an Edwards curve $E$ over a finite field $\mathbb{F}_p$ can be written in terms of the Gaussian hypergeometric series  $\displaystyle {_2F_1}\left(\begin{matrix}
                          \phi&\phi \\
                          {} & \epsilon
                        \end{matrix}\Big| x\right)$. Consequently, we evaluate $|E(\Fp)|$ for some $E$. Then we prove that every Edwards curve is isogenous to a Legendre curve over $\Fp$.
 Finally, it turns out that the number of rational points on a twisted Edwards curve $E$ is described using special values of the hypergeometric series $\displaystyle {_2F_1}\left(\begin{matrix}
                          \phi&\phi \\
                          {} & \epsilon
                        \end{matrix}\Big| x\right)$.
This sets the stage for evaluating $|E(\Fp)|$ for some twisted Edwards curves $E$.

  \section{Gaussian hypergeometric series}
 Throughout the note $p$ will be an odd prime unless otherwise stated. We extend multiplicative characters $\chi$ defined over $\mathbb{F}_p^{\times}$ to $\mathbb{F}_p$ by setting $\chi(0)=0$. We write $\overline{\chi}$ to denote $1/\chi$. The trivial and quadratic characters will be denoted by $\epsilon$ and $\phi$ respectively. Let $J(A,B)$ denote the Jacobi sum \[J(A,B)=\sum_{x\in\Fp}A(x)B(1-x)\] where $A$ and $B$ are characters over $\Fp$. Let $A_0,A_1,\ldots,A_n$ and $B_1,\ldots,B_n$ be characters defined over $\mathbb{F}_p$. The {\em Gaussian hypergeometric series} is \[_{n+1}F_n\left(\begin{matrix}
                          A_0&A_1&\ldots&A_n \\
                          {} & B_1&\ldots&B_n
                        \end{matrix}\Big| x\right):=\frac{p}{p-1}\sum_{\chi}{{A_0\chi}\choose{\chi}}{{A_1\chi}\choose{B_1\chi}}\ldots{{A_n\chi}\choose{B_n\chi}}\chi(x)\]
where the sum is over all characters over $\mathbb{F}_p$ and \[{A\choose B}:=\frac{B(-1)}{p}J(A,\overline{B})=\frac{B(-1)}{p}\sum_{x\in\mathbb{F}_p}A(x)\overline{B}(1-x).\]

The following properties of the symbol $\displaystyle{A\choose B}$ can be found in \cite{Greene}.

\begin{Lemma}
\label{lem1}
For any characters $A$ and $B$ over $\mathbb{F}_p$, one has:
\begin{myitemize}
\item[a)] $\displaystyle A(1+x)=\delta(x)+\frac{p}{p-1}\sum_{\chi}{A\choose \chi}\chi(x)$ where $\delta(x)=1$ if $x=0$ and $\delta(x)=1$ if $x\ne 0$;
\item[b)] $\displaystyle \overline{A}(1-x)=\delta(x)+\frac{p}{p-1}\sum_{\chi}{A\chi\choose \chi}\chi(x)$ where $\delta(x)=1$ if $x=0$ and $\delta(x)=0$ if $x\ne 0$;
\item[c)] $\displaystyle {A\choose B}={A\choose A\overline{B}}$;
\item[d)] $\displaystyle {A\choose B}={B\overline{A}\choose B}B(-1)$;
\item[e)] $\displaystyle {A\choose \epsilon}={A\choose A}=-\frac{1}{p}+\frac{p-1}{p}\delta(A)$ where $\delta(A)=1$ if $A=\epsilon$ and $\delta(A)=0$ otherwise;
\item[f)] $\displaystyle {B^2\chi^2\choose \chi}={\phi B \chi\choose\chi}{B\chi\choose B^2\chi}{\phi\choose\phi B}^{-1}B\chi(4).$
\end{myitemize}
\end{Lemma}

\section{Rational points on Edwards curves}

 Let $E$ be an elliptic curve over a field $k$ with $\Char k\ne 2$ defined by an Edwards model \[x^2+y^2=a^2(1+x^2y^2),\textrm{ where }a^5-a\ne 0.\] Such an elliptic curve will be called an {\em Edwards curve}. If $(x_1,y_1)$ and $(x_2,y_2)$ are two points on $E$, then these two points add up to \[x_3=\frac{1}{a}\cdot\frac{x_1y_2+x_2y_1}{1+x_1x_2y_1y_2},\;y_3=\frac{1}{a}\cdot\frac{y_1y_2-x_1x_2}{1-x_1x_2y_1y_2}.\]
 There are two points at infinity, namely if we homogenize the defining equation, we get \[x^2z^2+y^2z^2=a^2(z^4+x^2y^2)\] and the points at infinity are $(x:y:z)\in\{(1:0:0),(0:1:0)\}$.

 We will need the following lemma to count rational points on an Edwards curve.
 \begin{Lemma}
 \label{lem2}
 Let $A$ be a character on $\Fp$ and $a\in\Fp^{\times}$. The following identities hold:
 \begin{myitemize}
 \item[a)] $\displaystyle{A^2\choose A}=\left\{\begin{array}{ll}
{\phi A\choose A}A(4)  & \textrm{if } A\ne\epsilon \\
\frac{p-2}{p} & \textrm{if } A=\epsilon \end{array}\right.$
 \item[b)] $\displaystyle\sum_{x\in\Fp}A(a^2-x^2)=pA(4a^2){\overline{A}^2\choose\overline{A}}=\left\{\begin{array}{ll}
pA(a^2){\phi \overline{A}\choose \overline{A}}  & \textrm{if } A\ne\epsilon \\
p-2 & \textrm{if } A=\epsilon \end{array}\right.$
 \end{myitemize}
 \end{Lemma}
 \begin{Proof}
 a) In Lemma \ref{lem1} f), put $B=\epsilon$ and $\chi=A$. This yields \[{A^2\choose A}={\phi A\choose A}{A\choose A}{\phi\choose\phi}^{-1}A(4).\] According to Lemma \ref{lem1} e), the product above is $\displaystyle{\phi A\choose A}A(4)$ if $A\ne\epsilon$, and it is $\displaystyle(2-p){\phi A\choose A}A(4)$ if $A=\epsilon$.

 To prove b), we notice that
 \begin{eqnarray*}
 \sum_{x\in\Fp}A(a^2-x^2)&=&A\left(4a^2\right)\sum_{x\in\Fp}A\left(\frac{1}{2}-\frac{x}{2a}\right)A\left(\frac{1}{2}+\frac{x}{2a}\right)
 \end{eqnarray*}
 Setting $\displaystyle u=\frac{1}{2}-\frac{x}{2a}$, the above sum becomes
 \begin{eqnarray*}
  \sum_{x\in\Fp}A(a^2-x^2)&=&A\left(4a^2\right)\sum_{u\in\Fp}A(u)A(1-u)\\
   &=& A(4a^2)J(A,A)=pA(-4a^2){A\choose\overline{A}}.
 \end{eqnarray*}
 Using Lemma \ref{lem1} d), one has $\displaystyle {A\choose \overline{A}}={\overline{A}^2\choose\overline{A}}\overline{A}(-1)$. Part b) now follows from a).
 \end{Proof}
The following theorem relates the number of rational points on an Edwards curves over $\mathbb{F}_p$ to a Gaussian hypergeometric series.
\begin{Theorem}
\label{thm1}
Let $E/\Fp$ be described by $x^2+y^2=a^2(1+x^2y^2)$ where $a^5\not\equiv a$ mod $p$, $p$ is an odd prime. Then
\[|E(\mathbb{F}_p)|=1+p+p\cdot{_2F_1}\left(\begin{matrix}
                          \phi&\phi \\
                          {} & \epsilon
                        \end{matrix}\Big| 1-a^4\right).\]
\end{Theorem}
\begin{Proof}
 The defining equation of $E$ can be written as:
\[y^2=\frac{a^2-x^2}{1-a^2x^2}.\]
Bearing in mind that there are two points at infinity, one has
\begin{eqnarray*}
|E(\mathbb{F}_p)|&=&2+p+\sum_{x\in\Fp\setminus\{\pm a^{-1}\}}\phi\left(\frac{a^2-x^2}{1-a^2x^2}\right)\\
&=&2+p+\sum_{x\in\Fp\setminus\{\pm a^{-1}\}}\phi\left(1+\frac{a^2-a^{-2}}{a^{-2}-x^2}\right)\\
&=&2+p+\sum_{x\in\Fp\setminus\{\pm a^{-1}\}}\Big[\delta\left(\frac{a^2-a^{-2}}{a^{-2}-x^2}\right)+\frac{p}{p-1}\sum_{\chi}{\phi\choose\chi}\chi\left(\frac{a^2-a^{-2}}{a^{-2}-x^2}\right)\Big]\\
&=&2+p+\frac{p}{p-1}\sum_{x\in\Fp\setminus\{\pm a^{-1}\}}\sum_{\chi}{\phi\choose\chi}\chi\left(\frac{a^2-a^{-2}}{a^{-2}-x^2}\right)\\
&=&2+p+\frac{p}{p-1}\sum_{\chi}{\phi\choose\chi}\chi\left(a^2-a^{-2}\right)\sum_{x\in\Fp}\overline{\chi}(a^{-2}-x^2).\\
\end{eqnarray*}
The third equality follows from Lemma \ref{lem1} a). Now we use Lemma \ref{lem2} b) to obtain the following identity:
\begin{eqnarray*}
|E(\mathbb{F}_p)|&=&2+p+\frac{p}{p-1}\Big[p\sum_{\chi\ne\epsilon}{\phi\choose\chi}{\phi\chi\choose\chi}\overline{\chi}\left(a^{-2}\right)\chi\left(a^2-a^{-2}\right)+{\phi\choose\epsilon}(p-2)\Big].
\end{eqnarray*}
Lemma \ref{lem1} d) gives $\displaystyle{\phi\choose \chi}={\overline{\phi}\chi\choose \chi}\chi(-1)={\phi\chi\choose \chi}\chi(-1)$, thus
\begin{eqnarray*}
|E(\Fp)|&=&2+p+\frac{p^2}{p-1}\Big[\sum_{\chi\ne\epsilon}{\phi\chi\choose\chi}{\phi\chi\choose\chi}\chi\left(1-a^4\right)+\frac{2-p}{p^2}\Big]\\
&=&2+p+\frac{p^2}{p-1}\,\Big[\sum_{\chi}{\phi\chi\choose\chi}{\phi\chi\choose\chi}\chi\left(1-a^4\right)-{\phi\choose\epsilon}^2+\frac{2-p}{p^2}\Big]\\
&=&1+p+p\cdot{_2F_1}\left(\begin{matrix}
                          \phi&\phi \\
                          {} & \epsilon
                        \end{matrix}\Big|1-a^4\right).
\end{eqnarray*}
\end{Proof}
The Legendre family of elliptic curves is the one defined by \[E_{\lambda}:y^2=x(x-1)(x-\lambda),\;\lambda(\lambda-1)\not\equiv 0\textrm{ mod }p.\] Theorem 1 in \cite{OnoGauss} states that $\displaystyle |E_{\lambda}(\Fp)|=1+p+p\,\phi(-1)\cdot{_2F_1}\left(\begin{matrix}
                          \phi&\phi \\
                          {} & \epsilon
                        \end{matrix}\Big| \lambda\right)$. Since two elliptic curves over $\Fp$ are isogenous if and only if they have the same number of rational points, comparing $|E_{\lambda}(\Fp)|$ to the number of rational points of an Edwards curve $E$ over $\Fp$ defined by $x^2+y^2=a^2(1+x^2y^2)$, see Theorem \ref{thm1}, yields that $E$ is isogenous to the Legendre elliptic curve $E_{1-a^4}$ if $p\equiv 1$ mod $4$. In fact, every Edwards curve is isogenous to a Legendre curve over $\Fp$.

\begin{Corollary}
\label{cor:iso}
Let $E$ be defined by $x^2+y^2=a^2(1+x^2y^2)$ over $\Fp$, where $p$ is an odd prime and $a^5\not\equiv a$ mod $p$. Then $E$ is isogenous to the Legendre elliptic curve $E_{a^4}:y^2=x(x-1)(x-a^4)$.
\end{Corollary}
\begin{Proof}
According to Theorem 1 in \cite{OnoGauss}, $\displaystyle|E_{a^4}(\Fp)|=1+p+p\phi(-1)\cdot{_2F_1}\left(\begin{matrix}
                          \phi&\phi \\
                          {} & \epsilon
                        \end{matrix}\Big| a^4\right)$. The following identity is Theorem 4.4 (i) of \cite{Greene}
                        \[{_2F_1}\left(\begin{matrix}
                          \phi&\phi \\
                          {} & \epsilon
                        \end{matrix}\Big| x\right)=\phi(-1)\cdot{_2F_1}\left(\begin{matrix}
                          \phi&\phi \\
                          {} & \epsilon
                        \end{matrix}\Big| 1-x\right).\]
                       Now set $x=a^4$ and use Theorem \ref{thm1}. Consequently, $E$ and $E_{a^4}$ have the same number of rational points over $\Fp$. It follows that they are isogenous over $\Fp$.
\end{Proof}
We recall the following Proposition which can be found as Theorem 2 in \cite{OnoGauss}.
\begin{Proposition}
\label{prop1}
Let $p$ be an odd prime. If $\displaystyle\lambda\in\left\{-1,1/2,2\right\}$, then \begin{align*}{_2F_1}\left(\begin{matrix}
                          \phi&\phi \\
                          {} & \epsilon
                        \end{matrix}\Big| \lambda\right)=\left\{\begin{array}{ll}
0  & \textrm{if  } p\equiv 3 \textrm{ mod }4, \\
\frac{2x\cdot(-1)^{(x+y+1)/2}}{p} & \textrm{if } p\equiv 1\textrm{ mod }4,x^2+y^2=p, x\textrm{ is odd}. \end{array}\right.\end{align*}
\end{Proposition}
\begin{Corollary}
Let $E/\Fp$ be defined by $x^2+y^2=a^2(1+x^2y^2)$ where $a^5\not\equiv a$ mod $p$. If $a^4\in\{-1,1/2,2\}$, then \begin{align*}|E(\Fp)|=\left\{\begin{array}{ll}
1+p  & \textrm{if  } p\equiv 3 \textrm{ mod }4, \\
1+p+2x\cdot(-1)^{(x+y+1)/2} & \textrm{if } p\equiv 1\textrm{ mod }4,x^2+y^2=p, x\textrm{ is odd}. \end{array}\right.\end{align*}
\end{Corollary}
\begin{Proof}
This follows immediately from Corollary \ref{cor:iso}.
\end{Proof}

\section{Rational points on twisted Edwards curves}

Twisted Edwards curves were introduced in \cite{Bernstein} as generalizations of Edwards curves. These curves include more elliptic curves over finite fields than Edwards curves do. A twisted Edwards curve is an elliptic curve defined by the following {\em twisted Edwards model}
\[ax^2+y^2=1+dx^2y^2\] where $ad(a-d)\ne 0$. For any field $k$ with $\Char(k)\ne 2$, any elliptic curve over $k$ with three $k$-rational points of order $2$ is $2$-isogenous over $k$ to a twisted Edwards curve, see Theorem 5.1 of \cite{Bernstein}, hence they have the same number of rational points over $\Fp$.

\begin{Lemma}
\label{lem3}
The following equality holds for any character $A$ over $\Fp$
\[\sum_{x\in\Fp}\chi(x^2)\phi(1-x^2)=p\phi(-1)\left[{\phi\chi\choose\chi}+{\chi\choose\phi\chi}\right].\]
\end{Lemma}
\begin{Proof}
We recall that the number of solutions of the equation $z^2=a$ mod $p$ is given by $N_a=\phi(a)+1$. In fact $N_a=2$ if $a\in\Fp^2,a\ne0,$ $N_0=1$, and $N_a=0$ otherwise. We consider the following sum:
\begin{eqnarray*}
\sum_x\chi(x)\phi(x)\phi(1-x)&=&\sum_x\chi(x)\phi(1-x)\left(N_x-1\right)\\
&=&\sum_x\chi(x)\phi(1-x)N_x-\sum_x\chi(x)\phi(1-x)\\
&=&2\sum_{x\in\mathbb{F}_p^2}\chi(x)\phi(1-x)-\sum_x\chi(x)\phi(1-x)\\
&=&\sum_{x}\chi(x^2)\phi(1-x^2)-\sum_x\chi(x)\phi(1-x).\\
\end{eqnarray*}
Therefore,
\begin{eqnarray*}
\sum_x\chi(x^2)\phi(1-x^2)&=&J(\phi\chi,\phi)+J(\chi,\phi)\\
&=& p\phi(-1)\left[{\phi\chi\choose\phi}+{\chi\choose\phi}\right].
\end{eqnarray*}
Using Lemma \ref{lem1} c), one has $\displaystyle {\phi\chi\choose \phi}={\phi\chi\choose\chi}$, similarly $\displaystyle{\chi\choose\phi}={\chi\choose\phi\chi}$.
\end{Proof}
\begin{Theorem}
\label{thm2}
Let $E/\Fp$ be described by $ax^2+y^2=1+dx^2y^2$ where $ad(a-d)\not\equiv 0$ mod $ p$, $p$ is an odd prime. Then
\[|E(\mathbb{F}_p)|=2+p+\phi(a)+p\phi(-a)\left[{_2F_1}\left(\begin{matrix}
                          \phi&\phi \\
                          {} & \epsilon
                        \end{matrix}\Big| a^{-1}d\right)+{_2F_1}\left(\begin{matrix}
                          \phi&\epsilon \\
                          {} & \phi
                        \end{matrix}\Big| a^{-1}d\right)\right].\]
In particular, \[|E(\Fp)|=2+p-\phi(d)+p\phi(-a)\cdot{_2F_1}\left(\begin{matrix}
                          \phi&\phi \\
                          {} & \epsilon
                        \end{matrix}\Big| a^{-1}d\right).\]

\end{Theorem}
\begin{Proof}
We observe that we can write the above twisted Edwards model as \[x^2=\frac{1-y^2}{a-dy^2}.\]
We set $b=a^{-1}d$. The number of rational points of $E$ over $\Fp$ is given by
\begin{eqnarray*}
|E(\mathbb{F}_p)|&=&2+p+\sum_{y\in\Fp\setminus\{\sqrt{ad^{-1}}\}}\phi\left(\frac{1-y^2}{a-dy^2}\right)=2+p+\sum_{y\in\Fp}\phi\left(1-y^2\right)\phi\left(a-dy^2\right)\\
&=&2+p+\phi(a)\sum_{y\in{\Fp}}\phi(1-y^2)\phi(1-by^2)\\
&=&2+p+\phi(a)\sum_{y\in{\Fp}}\phi(1-y^2)\left[\delta(by^2)+\frac{p}{p-1}\sum_{\chi}{\phi\chi\choose\chi}\chi(by^2)\right].
\end{eqnarray*}
The last equality follows from Lemma \ref{lem1} b). According to Lemma \ref{lem3}, one obtains
\begin{eqnarray*}
|E(\Fp)|&=&2+p+\phi(a)+\frac{p\phi(a)}{p-1}\sum_{\chi}{\phi\chi\choose\chi}\chi(b)\sum_{y\in{\mathbb{F}_p}}\chi(y^2)\phi(1-y^2)\\
&=&2+p+\phi(a)+\frac{p^2\phi(-a)}{p-1}\sum_{\chi}{\phi\chi\choose\chi}\left[{\phi\chi\choose\chi}+{\chi\choose\phi\chi}\right]\chi(b)\\
&=&2+p+\phi(a)+p\phi(-a)\left[{_2F_1}\left(\begin{matrix}
                          \phi&\phi \\
                          {} & \epsilon
                        \end{matrix}\Big| a^{-1}d\right)+{_2F_1}\left(\begin{matrix}
                          \phi&\epsilon \\
                          {} & \phi
                        \end{matrix}\Big| a^{-1}d\right)\right].
\end{eqnarray*}
However, Corollary 3.16 in \cite{Greene} indicates that
\begin{eqnarray*}{_2F_1}\left(\begin{matrix}
                          \phi&\epsilon \\
                          {} & \phi
                        \end{matrix}\Big| a^{-1}d\right)&=&{\phi\choose \phi}\phi (-a^{-1}d)\epsilon(1-a^{-1}d)-\frac{1}{p}\phi(-1)\epsilon(a^{-1}d)+\frac{p-1}{p}\phi(-1)\delta(1-a^{-1}d)\delta(\epsilon)\\
                        &=&-\frac{1}{p}\phi(-a^{-1}d)-\frac{1}{p}\phi(-1).
                        \end{eqnarray*}
                        This proves the theorem.
\end{Proof}

\begin{Corollary}
Let $p$ be an odd prime. Let $E/\Fp$ be defined by $ax^2+y^2=1+dx^2y^2$ where $ad(a-d)\not\equiv 0$ mod $p$. If $\lambda:=a^{-1}d\in\{-1,1/2,2\}$, then
\begin{align*}|E(\Fp)|=\left\{\begin{array}{ll}
2+p-\phi(d)& \textrm{if  } p\equiv 3 \textrm{ mod }4, \\
2+p-\phi(d)+2x\cdot\phi(a)\cdot(-1)^{(x+y+1)/2} & \textrm{if } p\equiv 1\textrm{ mod }4,x^2+y^2=p, x\textrm{ is odd}. \end{array}\right.\end{align*}
\end{Corollary}
\begin{Proof}
This follows from Theorem \ref{thm2} and Proposition \ref{prop1}.
\end{Proof}

We remark that any Legendre curve is isogenous to a twisted Edwards curve over $\Fp$. Indeed, any elliptic curve $E$ with three $\Fp$-rational points of order $2$ defined by $y^2=x(x-a)(x-b)$ is isogenous to the twisted Edwards curve $4ax^2+y^2=1+4bx^2y^2$, see Theorem 5.1 in \cite{Bernstein}. Thus the formula for the number of rational points on a Legendre elliptic curve over $\Fp$, Theorem 1 of \cite{OnoGauss}, follows as a special case from Theorem \ref{thm2}.
\bibliographystyle{plain}
\footnotesize
\bibliography{Edwards}
Department of Mathematics and Actuarial Science\\ American University in Cairo\\ mmsadek@aucegypt.edu\\nelsissi@aucegypt.edu
\end{document}